\documentclass[a4paper]{amsart}

\usepackage[all]{xy}
\usepackage{amsmath,amsthm,amssymb,mathrsfs}
\usepackage{enumitem}
\usepackage[%
colorlinks=true,%
citecolor=green!70!black,%
linkcolor=red!80!black%
]{hyperref}
\usepackage{xcolor}
\usepackage{tikz}
\usetikzlibrary{arrows.meta}

\newtheorem{thm}{Theorem}[section]
\newtheorem{prop}[thm]{Proposition}
\newtheorem{cor}[thm]{Corollary}
\newtheorem{lem}[thm]{Lemma}

\theoremstyle{definition}
\newtheorem{ex}[thm]{Example}
\newtheorem{defn}[thm]{Definition}
\newtheorem{rem}[thm]{Remark}
\newtheorem*{notation}{Notation}

\renewcommand{\div}{\operatorname{div}}
\newcommand{\del}{\partial}

\title{Approximation and extension theorems for normal affine toric varieties}
\author{Ryo Uchimoto}
\date{\today}
\address{
Joint Graduate School of Mathematics for Innovation\\
Kyushu University \\
744 Motooka \\
Nishi-ku \\
Fukuoka 819-0395 \\
Japan 
}
\email{uchimoto.ryo.387@s.kyushu-u.ac.jp}


\begin{document}


\begin{abstract}
We study approximation and extension problems for holomorphic maps from
Stein spaces to normal affine toric varieties, continuing the work of
L\"ark\"ang and L\'arusson. Normality allows us to obtain twisted
factorizations for arbitrary holomorphic maps from factorial complex
spaces, including degenerate maps. This yields approximation and extension
theorems under weaker cohomological assumptions. We also prove a simultaneous approximation and extension theorem
under an additional codimension condition on the inverse image of the
singular locus and a cohomological condition on the union of the
interpolation submanifold and the approximation set. Finally, an explicit
example with target $\{w^2=uv\}\cong\mathbb{C}^2/\{\pm1\}$ shows that simultaneous
approximation and extension can fail even when the relevant cohomological
obstructions vanish. Thus normality alone does not recover the full
approximation and interpolation property of Oka manifolds.
\end{abstract}

\maketitle

\tableofcontents 


\section{Introduction}\label{sec:intro}

Oka theory studies approximation and extension problems for holomorphic
maps from Stein spaces to suitable complex analytic targets, with the aim of
identifying situations in which topological obstructions are the only
obstructions to a holomorphic solution. Its origins include Oka's theorem
that a topologically trivial holomorphic line bundle on a Stein manifold
is holomorphically trivial~\cite{Oka39}. Grauert generalized this principle
to principal bundles with complex Lie structure group~\cite{Gra58}, and
Gromov developed it further through the theory of elliptic complex
manifolds and holomorphic sprays~\cite{Gro89}.

Two central properties in modern Oka theory are the \emph{Convex
Approximation Property} (CAP) and the \emph{Parametric Oka Property with
Approximation and Interpolation} (POPAI). A complex manifold~$Y$ has CAP
if, for every $n\geqslant1$ and every compact convex set $K\subset\mathbb{C}^n$, every
holomorphic map from a neighborhood of~$K$ to~$Y$ can be approximated
uniformly on~$K$ by entire maps $\mathbb{C}^n\to Y$.

POPAI combines this kind of approximation with extension from a closed
complex subvariety, allowing continuous parameters and a homotopy. More
precisely, let $X$ be a reduced Stein space, let $X'\subset X$ be a closed
complex subvariety, and let $K\subset X$ be a compact $\mathcal{O}(X)$-convex set.
For compact parameter spaces $P_0\subset P\subset\mathbb R^m$, consider a
continuous map
\[
 f\colon P\times X\to Y
\]
such that $f(p,\cdot)$ is holomorphic on~$X'$ and on a common neighborhood
of~$K$ for every $p\in P$, and is holomorphic on all of~$X$ for
$p\in P_0$. POPAI asserts that $f$ can be deformed through a continuous
homotopy $f_t\colon P\times X\to Y$, $t\in[0,1]$, with $f_0=f$, such that
\begin{enumerate}[label=\textup{(\roman*)},leftmargin=*,itemsep=2pt]
\item $f_1(p,\cdot)$ is holomorphic on~$X$ for every $p\in P$;
\item each $f_t(p,\cdot)$ remains holomorphic near~$K$, and the homotopy
is uniformly as close to~$f$ on $P\times K$ as desired;
\item $f_t=f$ on $(P_0\times X)\cup(P\times X')$ for every $t\in[0,1]$.
\end{enumerate}
When $P$ is a point and $P_0=\varnothing$, this is the \emph{Basic Oka
Property with Approximation and Interpolation} (BOPAI).

Extending a given holomorphic map on~$X'$ means constructing a holomorphic
map on~$X$ that agrees with it on~$X'$. This agreement condition is also
referred to as \emph{interpolation along~$X'$}. In the simultaneous
problem, we seek such an extension that also approximates given
holomorphic data near~$K$, uniformly on~$K$. We therefore describe this
problem both as \emph{simultaneous approximation and extension} and as
\emph{approximation with interpolation along~$X'$}.

Forstneri\v{c} proved that CAP implies POPAI~\cite{For09}. These
properties, together with their jet interpolation variants and several
other ostensibly different properties, are equivalent for complex
manifolds and characterize the class of \emph{Oka manifolds}; see
\cite[\S5.15]{For17}. Gromov's $\text{Ell}_\infty$ property
\cite[\S3.1]{Gro89} is the formulation with approximation and finite-order
jet interpolation for finite polyhedral parameter pairs. It assumes that
the initial maps are holomorphic on a neighborhood of the interpolation
subvariety as well as of the approximation set. For manifold targets it
is equivalent to POPAI; see \cite[\S5.15]{For17}. For singular targets,
however, holomorphicity on~$X'$ need not give a holomorphic extension to a
neighborhood of~$X'$, so these formulations must be distinguished.

\medskip

Developing an Oka theory for singular targets is a natural problem raised
by Gromov~\cite{Gro89}; see also~\cite{For17}. In particular, Gromov
suggested the quotient $\mathbb{C}^2/\{\pm1\}$, where $-1$ acts by
$(s,t)\mapsto(-s,-t)$, as a candidate singular $\text{Ell}_\infty$ space.
This quotient is the normal affine toric surface
\[
 \{(u,v,w)\in\mathbb{C}^3:w^2=uv\}.
\]
Affine toric varieties therefore provide a concrete class of singular
targets in which to investigate which approximation and extension
properties survive.

L\"ark\"ang and L\'arusson~\cite{LL16} initiated this investigation for
affine toric varieties that are not necessarily normal. Every affine
toric variety satisfies CAP, by the standard argument used in the proof
of \cite[Proposition~1.5]{LL16}. In contrast, they showed that for the
Whitney umbrella
\[
 Y=\{(x,y,z)\in\mathbb{C}^3:x^2y=z^2\}
\]
there is a holomorphic map $\mathbb{C}\times\{0\}\cup\{(0,1)\}\to Y$ with no
holomorphic extension to~$\mathbb{C}^2$, although it has a continuous extension.
Thus CAP does not imply POPAI for singular targets. Their negative results
also include nonextendable maps from smooth surfaces biholomorphic to a
product of two annuli into singular normal affine toric varieties without
torus factors~\cite[Theorem~1.4\textup{(b)}]{LL16}.
On the positive side, they proved the following extension theorem.
Recall that a complex space is \emph{factorial} if all of its local rings
are unique factorization domains; in particular, every complex manifold
is factorial.

\begin{thm}\label{thm:LL}
Let $Y$ be a not necessarily normal affine toric variety. Let $X'$ be a closed
factorial subvariety of a reduced Stein space $X$ such that
$H^q(X,\mathbb{Z})\to H^q(X',\mathbb{Z})$ is surjective for $q=0,1,2$.
Then every holomorphic map $X'\to Y$ extends to a holomorphic map $X\to Y$.
\end{thm}

The first aim of this paper is to clarify the role of normality in these
results. The main tool of~\cite{LL16} is the \emph{twisted factorization}
of a nondegenerate holomorphic map into an affine toric variety, expressed
in terms of sections of holomorphic line bundles. Here a map is called
\emph{nondegenerate} if the image of each irreducible component of its
source meets the dense torus of the target. For a normal affine toric
variety, we show that a twisted factorization can be constructed even for
degenerate maps, using fixed factorization data for the ambient target
(Theorem~\ref{thm:tfact}). This gives the following extension theorem under weaker
cohomological assumptions than those in Theorem~\ref{thm:LL}.

\begin{thm}\label{thm:extension}
Let $Y$ be a normal affine toric variety. Let $X'$ be a closed factorial
subvariety of a reduced Stein space $X$ such that
$H^q(X,\mathbb{Z})\to H^q(X',\mathbb{Z})$ is surjective for $q=1,2$.
Then every holomorphic map $X'\to Y$ extends to a holomorphic map $X\to Y$.
\end{thm}

The same factorization method also yields an approximation theorem.

\begin{thm}\label{thm:approximation}
Let $Y$ be a normal affine toric variety. Let $\Omega$ be a Stein Runge open
subset of a factorial Stein space $X$ such that
$H^q(X,\mathbb{Z})\to H^q(\Omega,\mathbb{Z})$ is surjective for $q=1,2$.
Then every holomorphic map $\Omega\to Y$ can be approximated uniformly on
compact subsets of $\Omega$ by holomorphic maps $X\to Y$.
\end{thm}

The cohomological assumptions in these theorems ensure that the line
bundles arising in twisted factorizations extend to~$X$ and that the
torus factors of~$Y$ introduce no topological obstruction. Although approximation and extension
can be obtained separately, their simultaneous realization requires
compatible factorizations of the data near the approximation set and on
the interpolation subvariety. The second aim of this paper is to identify
a condition ensuring this compatibility.

We show that the twisted factorization of a nondegenerate holomorphic map
from a complex manifold is unique up to isomorphism if the inverse image
of the singular locus has codimension at least~$2$ (Lemma~\ref{lem:comp}). Applying
this on the overlap of the two sets gives the following simultaneous
approximation and extension theorem.

\begin{thm}\label{thm:simultaneous}
Let $X$ be a Stein manifold, $X'\subset X$ be a closed submanifold, and $K\subset X$
be an $\mathcal{O}(X)$-convex compact set. Let $Y$ be a normal affine toric variety.
Let $f\colon X\to Y$ be a continuous map which is holomorphic on $X'$ and
in some neighborhood $U$ of $K$. Assume that
\begin{enumerate}[label={\rm (\arabic*)}]
\item $f$ is nondegenerate on $X'\cap U$ and
$\operatorname{codim}_{X'\cap U}f^{-1}(Y_{\mathrm{sing}})\geqslant 2;$
\item the restriction map $H^2(X,\mathbb{Z})\to\check H^2(X'\cup K,\mathbb{Z})$ is surjective.
\end{enumerate}
Then there is a sequence of holomorphic maps $F_j\colon X\to Y$, $j=1,2,\ldots$
such that $F_j|_{X'}=f|_{X'}$ and $F_j\to f$ uniformly on $K$ as $j\to\infty$.
Moreover, for each $j$ there is a homotopy
$f_{j,\theta}\colon X\to Y$ $(0\leqslant\theta\leqslant1)$ with
$f_{j,0}=f$ and $f_{j,1}=F_j$ such that, for every $\theta$,
$f_{j,\theta}|_{X'}=f|_{X'}$ and $f_{j,\theta}$ is holomorphic on a
neighborhood of $K$, and $f_{j,\theta}\to f$ as $j\to\infty$ uniformly
on $K$ and uniformly in $\theta\in[0,1]$.
\end{thm}

Here $\check H^2$ denotes \v{C}ech cohomology, whose continuity property
is used in the proof. Singular cohomology does not satisfy this property
in general, and $X'\cup K$ need not be locally contractible, so the two
cohomology theories cannot be identified a priori.
The homotopies in Theorem~\ref{thm:simultaneous} preserve the
interpolation data and remain holomorphic near~$K$, with approximation
uniform in the homotopy parameter.

Finally, Example~\ref{ex:notAI} shows that the codimension condition in
Theorem~\ref{thm:simultaneous} cannot simply be omitted. For
$X=\mathbb{C}^3$, $X'=\mathbb{C}^2\times\{0\}$, a suitable polynomially convex
polydisc~$K$, and
$Y=\{w^2=uv\}\cong\mathbb{C}^2/\{\pm1\}$, we construct a continuous map
$f\colon X\to Y$, holomorphic on~$X'$ and near~$K$, for which no sequence
of holomorphic maps $X\to Y$ can both agree with~$f$ on~$X'$ and converge
to~$f$ uniformly on~$K$. In this example the relevant cohomological
obstructions vanish, and the prescribed map $f|_{X'}$ even has a global
holomorphic extension. The obstruction concerns compatibility with the
approximation data. Thus even for normal affine toric targets, the
simultaneous approximation and interpolation property can fail in a
setting where the separate problems are solvable.

The paper is organized as follows. Section~\ref{sec:toric} recalls the required facts
about affine toric varieties. Section~\ref{sec:factorization} develops twisted factorization
for normal targets and proves Theorems~\ref{thm:extension}
and~\ref{thm:approximation}. Section~\ref{sec:BOPAI} proves the uniqueness lemma and
Theorem~\ref{thm:simultaneous}, and concludes with the counterexample.


\section{Basic facts from the theory of affine toric varieties}\label{sec:toric}

In this paper we always assume that $\mathbb{N}=\mathbb{Z}_{\geqslant0}$.

\subsection{Definition of affine toric varieties}\label{subsec:affine}

Let $N$ be a free $\mathbb Z$-module of rank $d$, and put $M=\operatorname{Hom}(N,\mathbb Z)$.
We write $N_{\mathbb R}=N\otimes_{\mathbb Z}\mathbb R$ and
$M_{\mathbb R}=M\otimes_{\mathbb Z}\mathbb R$.
For an affine semigroup $\mathsf S\subset M$ with $\mathbb Z\mathsf S=M$
(see \cite[p.~16]{CLS11} for the definition), we associate the 
$\mathbb C$-algebra 
\[\mathbb{C}[\mathsf{S}]=\left\{\sum_{u\in\mathsf{S}}c_u\chi^u\;\middle|\;
c_u\in\mathbb{C}~\text{and}~c_u=0~\text{for all but finitely many $u$}\right\}\]
with multiplication induced by $\chi^u\cdot\chi^{u'}=\chi^{u+u'}$. 
Since an affine semigroup $\mathsf{S}$ is finitely generated, 
$\mathbb{C}[\mathsf{S}]$ is a finitely generated $\mathbb{C}$-algebra. 
Now let $Y_{\mathsf{S}}:=\operatorname{Spec}(\mathbb{C}[\mathsf{S}])^{\rm an}$. 
We call $Y_\mathsf{S}$ the \emph{affine toric variety} 
associated with the affine semigroup $\mathsf{S}$.\\

The normality of $Y_\mathsf{S}$ can be determined 
by whether $\mathsf{S}$ is saturated.

\begin{defn}
We say an affine semigroup $\mathsf{S}$ is saturated if 
for any $a\in\mathbb{Z}_{>0}$ and $u\in M$, $au\in\mathsf{S}$ implies $u\in\mathsf{S}$.
\end{defn}
\begin{defn}
\begin{enumerate}[label=(\arabic*)]
\item A \emph{rational polyhedral cone} $\sigma\subset N_\mathbb{R}$ 
is a set of the form
\[\sigma=\left\{\sum_{v\in A}\lambda_v v\;\middle|\;\lambda_v\geqslant 0\right\}\]
for some finite subset $A\subset N$.
\item A rational polyhedral cone $\sigma$ is \emph{strongly convex} 
if $\sigma\cap(-\sigma)=\{0\}$.
\item A rational polyhedral cone $\sigma$ is \emph{full dimensional} 
if $\dim\sigma=\dim N_{\mathbb R}=d$.
\end{enumerate}
\end{defn}

\begin{thm}[{\cite[Theorem 1.3.5]{CLS11}}]
Let $Y_\mathsf{S}$ be an affine toric variety associated 
with the affine semigroup $\mathsf{S}$. Then the following are equivalent.
\begin{enumerate}[label={\rm (\alph*)}]
\item $Y_\mathsf{S}$ is normal.
\item $\mathsf{S}$ is saturated.
\item There exists a strongly convex rational polyhedral cone 
$\sigma\subset N_\mathbb{R}$ such that $\mathsf{S}=\sigma^\lor\cap M$, 
where $\sigma^\lor=\{u\in M_\mathbb{R}\mid \langle u,v\rangle\geqslant0,~\text{\rm for~all}~v\in\sigma\}$.
\end{enumerate}
\end{thm}

For a strongly convex rational polyhedral cone $\sigma\subset N_{\mathbb R}$, we write
$\mathsf S_\sigma=\sigma^\lor\cap M$ and 
$Y_\sigma=\operatorname{Spec}(\mathbb C[\mathsf S_\sigma])^{\rm an}$.
It is well known that $\sigma\cap N$ has a unique minimal generating set.
This generating set is called the \emph{Hilbert basis} of $\sigma\cap N$. 
Note that $\sigma^\lor$ is strongly convex if and only if 
$\sigma$ is full dimensional (\cite[Proposition 1.2.12]{CLS11}).
In what follows for a full dimensional strongly convex rational 
polyhedral cone $\sigma\subset N_\mathbb{R}$ 
we always regard $Y_\sigma$ as embedded in $\mathbb{C}^n$
via the surjective $\mathbb{C}$-algebra homomorphism
\[
\mathbb{C}[x_1,\dots,x_n]\to\mathbb{C}[\sigma^\lor\cap M];\quad x_j\mapsto \chi^{u_j}
\]
where $\{u_1,\dots,u_n\}$ is the Hilbert basis of $\sigma^\lor\cap M$.

Geometrically, full-dimensionality can be characterized as follows.

\begin{prop}\label{prop:0inY}
Let $\sigma \subset N_{\mathbb{R}}$ be a strongly convex rational 
polyhedral cone, and put $\mathsf{S}_\sigma=\sigma^\lor\cap M$. 
Choose nonzero elements $u_1,\dots,u_n$ generating $\sigma^\lor\cap M$, 
and regard $Y_\sigma$ as embedded in  $\mathbb{C}^n$ via $x_j \mapsto \chi^{u_j}$. 
Then the following are equivalent:
\begin{enumerate}[label={\rm (\arabic*)}]
\item $\sigma$ is full-dimensional.
\item $Y_\sigma$ has no torus factor.
\item $Y_\sigma$ contains the origin of $\mathbb{C}^n$.
\end{enumerate}
\end{prop}

\begin{proof}
The equivalence of (1) and (2) follows from \cite[Proposition 3.3.9]{CLS11}.
Explicitly, set $N_\sigma=N\cap\operatorname{Span}_{\mathbb{R}}(\sigma)$. 
Choosing a complement of $N_\sigma$ in $N$ gives a toric isomorphism
\[
Y_\sigma \simeq Y_{\sigma,N_\sigma} \times (\mathbb{C}^*)^{\operatorname{rank}N-\dim\sigma},
\]
where $Y_{\sigma,N_\sigma}$ denotes the affine toric variety 
associated with $\sigma$ considered as a cone in $N_{\sigma,\mathbb{R}}$. 
Since $\sigma$ is full-dimensional in $N_{\sigma,\mathbb{R}}$, 
the first factor has no torus factor.

Suppose that (1) holds. Then $\sigma^\lor$ is strongly convex, 
so $\mathsf{S}_\sigma\cap(-\mathsf{S}_\sigma)=\{0\}$. Consequently, the assignment
\[
\chi^m \mapsto
\begin{cases}
1, & \text{if } m=0,\\
0, & \text{if } m\ne 0
\end{cases}
\]
extends to a $\mathbb{C}$-algebra homomorphism $\mathbb{C}[\mathsf{S}_\sigma]\to\mathbb{C}$. 
The corresponding point has all coordinates zero, proving (3).

Conversely, if $\sigma$ is not full-dimensional, choose $0\neq m\in\sigma^\perp\cap M$. 
Both $m$ and $-m$ belong to $\mathsf{S}_\sigma$, so we can write 
$m=\sum_j a_j u_j$ and $-m=\sum_j b_j u_j$ with $a_j,b_j\in\mathbb{Z}_{\ge 0}$. Thus
\[
\prod_{j=1}^n x_j^{a_j+b_j}=1
\]
on $Y_\sigma$. Since some $a_j+b_j$ is positive, this equation excludes the origin. 
Hence (3) implies (1).
\end{proof}

\subsection{Torus and torus invariant subvarieties}\label{subsec:torus}
Let $\mathsf S\subset M$ be an affine semigroup with $\mathbb Z\mathsf S=M$.
Its \emph{torus} is
\[
T_N=\operatorname{Spec}(\mathbb C[M])^{\rm an}\cong(\mathbb C^*)^d.
\]
The inclusion $\mathbb C[\mathsf S]\hookrightarrow\mathbb C[M]$ induces an open embedding
$T_N\hookrightarrow Y_{\mathsf S}$, and we identify $T_N$ with its image.

\begin{prop}
Let $Y$ be an affine toric variety. The torus $T_N\subset Y$ is a 
Zariski open subset such that the action of $T_N$ on itself 
(induced by the multiplication of $(\mathbb{C}^\ast)^d$) 
extends to an algebraic action of $T_N$ on $Y$.
\end{prop}

We omit its proof (see \cite[Propositions 1.1.8, 1.1.14]{CLS11}). 

\begin{defn}
Let $Y'$ be a subvariety of an affine toric variety $Y$. 
We say $Y'$ is $T_N$-invariant if $T_N\cdot Y'\subset Y'$.
\end{defn}

Now let $\sigma\subset N_{\mathbb R}$ be a strongly convex rational polyhedral cone.
For a face $\tau\preceq\sigma$, put
\[
N_\tau=N\cap\operatorname{Span}_{\mathbb R}(\tau),\qquad N(\tau)=N/N_\tau,
\]
and identify $M(\tau)=\operatorname{Hom}(N(\tau),\mathbb Z)$ with $M\cap\tau^\perp$.
Write $N(\tau)_{\mathbb R}=N(\tau)\otimes_{\mathbb Z}\mathbb R$, and denote the quotient
map $N\to N(\tau)$ and its real-linear extension by $q$.
We associate a point $\gamma_\tau\in Y_\sigma$ with the surjective $\mathbb C$-algebra
homomorphism defined, for $m\in\mathsf S_\sigma$, by
\[\mathbb C[\mathsf S_\sigma]\to\mathbb C;\qquad\chi^m\mapsto\begin{cases}
1, & m\in\tau^\perp\cap M, \\
0, & \text{otherwise}.
\end{cases}\]
We define $O(\tau):=T_N\cdot \gamma_\tau$ and $V(\tau)=\overline{O(\tau)}$. 
Let us recall the Orbit-Cone Correspondence.

\begin{thm}[{\cite[Theorem 3.2.6]{CLS11}}]\label{thm:orbit-cone}
Let $Y$ be a normal affine toric variety associated 
with an affine semigroup $\mathsf{S}_\sigma$. 
Then $\tau\mapsto V(\tau)$ is an inclusion-reversing bijection 
between the faces of $\sigma$ and the irreducible 
$T_N$-invariant closed analytic subsets of $Y$. 
Moreover $V(\tau)$ is the normal affine toric variety 
associated with $\mathsf{S}_\sigma\cap\tau^\perp=q(\sigma)^\lor\cap M(\tau)$, 
and $\dim V(\tau)=d-\dim\tau$.
\end{thm}

\subsection{Smooth and singular loci}

\begin{defn}
Let $\sigma$ be a strongly convex rational polyhedral cone.
\begin{itemize}
\item Let $\rho$ be an edge of $\sigma$. 
The semigroup $\rho\cap N$ is generated by a unique element $v_\rho\in \rho\cap N$. 
$v_\rho$ is said to be the \emph{ray generator} of $\rho$.
\item $\sigma$ is said to be \emph{smooth} if its ray generators 
$\{v_\rho\mid \rho~\text{is an edge of $\sigma$}\}$ can be 
extended to a $\mathbb{Z}$-basis of $N$. 
In particular $v_\rho$'s are linearly independent over $\mathbb{Z}$.
\end{itemize}
\end{defn}

\begin{prop}[{\cite[Proposition 11.1.2]{CLS11}}]\label{prop:smoothsing}
Let $Y_{\sigma}$ be a normal affine toric variety. Then
\[(Y_\sigma)_{\rm sing}=\bigcup\{V(\tau)\mid \tau\preceq\sigma:{\rm not~smooth}\},
\quad (Y_\sigma)_{\rm sm}=\bigcup\{U_\tau\mid\tau\preceq\sigma:{\rm smooth}\}.\]
Here $\mathsf S_\tau=\tau^\lor\cap M$, and
$U_\tau=\operatorname{Spec}(\mathbb C[\mathsf S_\tau])^{\rm an}$ 
is the corresponding toric open subset of $Y_\sigma$.
\end{prop}

\begin{ex}
Let $\sigma\subset\mathbb{R}^2$ be the strongly convex 
rational polyhedral cone shown in figure~\ref{fig:sigma}.
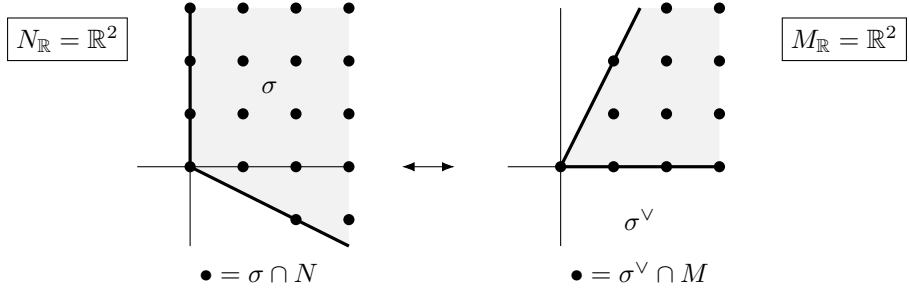
\begin{figure}[h]
\begin{tikzpicture}[scale=.7]

\fill[black!5] (0,3)--(0,0)--(3,-1.5)--(3,3)--(0,3);
\draw (-1,0)--(3,0);
\draw (0,-1.5)--(0,3); 
\draw[very thick] (0,3)--(0,0)--(3,-1.5);
\node at (1.5,1.5) {$\sigma$};
\foreach \i in {0,1,2,3}
\fill (0,\i) circle (3pt);
\foreach \i in {0,1,2,3}
\fill (1,\i) circle (3pt);
\foreach \i in {-1,0,1,2,3}
\fill (2,\i) circle (3pt);
\foreach \i in {-1,0,1,2,3}
\fill (3,\i) circle (3pt);
\node[below left] at (-1,3) {\fbox{$N_\mathbb{R}=\mathbb{R}^2$}};
\node[right] at (0,-2) {$\bullet=\sigma\cap N$};

\draw[Latex-Latex]  (4,0)--(5,0);

\fill[black!5] (7,0)--(10,0)--(10,3)--(8.5,3)--cycle;
\draw[very thick] (8.5,3)--(7,0)--(10,0);
\draw (6,0)--(10,0);
\draw (7,-1.5)--(7,3);
\foreach \i in {7,8,9,10}
\fill (\i,0) circle (3pt);
\foreach \i in {8,9,10}
\fill (\i,1) circle (3pt);
\foreach \i in {8,9,10}
\fill (\i,2) circle (3pt);
\foreach \i in {9,10}
\fill (\i,3) circle (3pt);
\node at (8.5,-1) {$\sigma^\lor$};
\node[below right] at (11,3) {\fbox{$M_\mathbb{R}=\mathbb{R}^2$}};
\node[right] at (7,-2) {$\bullet=\sigma^\lor\cap M$};
\end{tikzpicture}
\caption{Example of a strongly convex rational polyhedral cone}
\label{fig:sigma}
\end{figure}

It is easy to see that the Hilbert basis of $\mathsf{S}_\sigma=\sigma^\lor\cap M$ 
is $\{(1,0),(1,1),(1,2)\}$ so $\mathbb{C}[\mathsf{S}_\sigma]$ is 
generated by $\chi^{(1,0)},\chi^{(1,1)},\chi^{(1,2)}$. 
Let $\mathbb{C}[x,y,z]\to\mathbb{C}[\mathsf{S}_\sigma]$ be the surjective 
$\mathbb{C}$-algebra homomorphism defined by
\[x\mapsto \chi^{(1,0)},\quad y\mapsto \chi^{(1,2)},\quad z\mapsto\chi^{(1,1)}.\]
Since its kernel is $(z^2-xy)$, the corresponding normal affine toric variety 
$Y_\sigma$ is $\{z^2=xy\}$. 

There are four faces of $\sigma$, namely $\sigma$, $\mathbb{R}_{\geqslant0}(0,1)$, 
$\mathbb{R}_{\geqslant0}(2,-1)$, and $\{0\}$. We see that
\begin{align*}
&V(\sigma) = \{0\},\quad V(\{0\})=Y_\sigma, \\
&V(\mathbb{R}_{\geqslant0}(0,1)) = \mathbb{C}\times\{0\}^2,
\quad V(\mathbb{R}_{\geqslant0}(2,-1)) = \{0\}\times\mathbb{C}\times\{0\}.
\end{align*}
Note that the torus $T_N$ is $Y_\sigma\cap(\mathbb{C}^\ast)^3$.

The set of edges of $\sigma$ is $\{\mathbb{R}_{\geqslant0}(0,1),\mathbb{R}_{\geqslant0}(2,-1)\}$ 
so $\sigma$ is not smooth. Since the face $\mathbb{R}_{\geqslant0}(0,1)$ is itself 
a strongly convex rational polyhedral cone and its ray generator $(0,1)$ can be 
extended to a $\mathbb{Z}$-basis of $N$, it is smooth. The faces $\mathbb{R}_{\geqslant0}(2,-1)$ and $\{0\}$ 
are also smooth. Therefore we get
\[(Y_\sigma)_{\rm sing}=V(\sigma)=\{0\},\quad (Y_\sigma)_{\rm sm}=Y_\sigma\setminus\{0\}.\]
\end{ex}


\section{Approximation or extension theorems for normal affine toric varieties}\label{sec:factorization}

\subsection{Twisted factorization theorem}

We first recall the so-called \emph{twisted factorization} of a nondegenerate holomorphic map. 
Let $W$ be a factorial complex space and $Y\subset\mathbb{C}^n$ be a not necessarily normal affine toric variety. 
We use the toric embedding associated with a finite generating set of the defining semigroup.
We further assume that $0\in Y$ and the normalization of $Y$ is not isomorphic to affine space.

\begin{defn}
Let $X$ be a complex space and $Y$ be an affine toric variety. 
A holomorphic map $f\colon X\to Y$ is said to be nondegenerate
if for any irreducible component $X'\subset X$, 
the image $f(X')$ intersects the torus $T_N$ of $Y$.
\end{defn}

We first fix the matrices $A$, $B$, and $E$ associated with $Y$.

\paragraph{Construction of $A$}\label{constructionA}
Let $I$ be the toric ideal of $Y$. It is well known that there is a lattice $\Lambda\subset\mathbb{Z}^n$ such that
\[I=\langle x^{\alpha}-x^\beta
\mid \alpha,\beta\in\mathbb{N}^n~\text{and}~\alpha-\beta\in\Lambda\rangle.\]
Since $\mathbb{C}[x_1,\dots,x_n]$ is noetherian there are finitely many binomials $x^{\alpha_\nu}-x^{\beta_\nu}~(\nu=1,\dots,k)$ which generate $I$. 
We regard $\alpha_\nu,\beta_\nu$ as row vectors and put
\[
A=\begin{pmatrix}
\alpha_1-\beta_1 \\
\vdots \\
\alpha_k-\beta_k
\end{pmatrix}\in M_{k,n}(\mathbb{Z}).
\]

\paragraph{Construction of $B$}
The assumption that $0\in Y$ implies that $K=\operatorname{Ker}A\cap\mathbb{N}^n$ is nontrivial. Let $\mathscr{B}=\{b_1,\dots,b_m\}$ be the Hilbert basis of $K$, 
i.e., $\mathscr{B}$ is the smallest subset of $K$ with the following properties: $\mathscr{B}$ is a finite subset of $K$ and every vector in $K$ 
can be written as a (not necessarily unique) linear combination of the vectors in $\mathscr{B}$ with coefficients in $\mathbb{N}$. 
We regard the vectors in 
$\mathscr{B}$ as column vectors and put
\[
B=\begin{pmatrix}
b_1 & \cdots & b_m
\end{pmatrix}\in M_{n,m}(\mathbb{N}).
\]

\paragraph{Construction of $E$}
The assumption that the normalization of $Y$ is not $\mathbb{C}^d$ implies that $\operatorname{Ker}B\subset\mathbb{Z}^m$ is nontrivial. 
Since $\operatorname{Ker}B$ is free, it has a $\mathbb{Z}$-basis $\{e_1,\dots,e_l\}$. We naturally regard the vectors $e_\lambda~(\lambda=1,\dots,l)$ as column vectors and put
\[E=\begin{pmatrix}
e_1 & \cdots & e_l
\end{pmatrix}\in M_{m,l}(\mathbb{Z}).\]

By construction, $AB=0$ and $BE=0$. We keep these matrices fixed throughout this subsection.

\begin{thm}[{\cite[Theorem 2.2]{LL16}}]\label{thm:TFT}
Let $W,Y$ be as above, and let $A,B,E$ be the fixed matrices just constructed.
Then every nondegenerate holomorphic map $f\colon W\to Y$ has a twisted factorization:
there are holomorphic line bundles $N_\lambda$ $(\lambda=1,\dots,l)$ such that, for
$L_i=N_1^{e_{i1}}\otimes\dots\otimes N_l^{e_{il}}$ $(i=1,\dots,m)$,
there are global sections $\sigma_i\in H^0(W,L_i)$ satisfying
\[
f_j=\sigma_1^{b_{j1}}\otimes\dots\otimes\sigma_m^{b_{jm}},\qquad j=1,\dots,n,
\]
in the canonical trivializations.
\end{thm}

\begin{notation}
In the setting of Theorem~\ref{thm:TFT} we will write $L=N^E$ and $f=\sigma^B$. If we were to use additive notation for $\operatorname{Pic}(W)$ we would write $L=EN$.
\end{notation}

\subsection{Advantages of normality}\label{subsec:adv}
Next we further assume that $Y$ is normal, i.e., $Y$ arises from a full dimensional strongly convex rational polyhedral cone $\sigma\subset N_\mathbb{R}$.
We fix our notation as follows:
\begin{itemize}
\item $\tau\preceq\sigma$ is a face of $\sigma$; the lattices $N_\tau$, $N(\tau)$,
$M(\tau)$ and the projection $q$ are as in \S~\ref{subsec:torus}.
\item $\{v_i\mid i=1,\dots,m\}$ is the Hilbert basis of $\sigma\cap N$.
\item $\{v_i'\mid i=1,\dots,m'\}$ is the Hilbert basis of $q(\sigma)\cap N(\tau)=q(\sigma\cap N)$.
\item $\{u_j\mid j=1,\dots,n\}$ is the Hilbert basis of $\sigma^\lor\cap M$.
\item After reordering the $u_j$, the subset $\{u_1,\dots,u_{n'}\}$ is the Hilbert basis
of the face $\sigma^\lor\cap\tau^\perp\cap M$ of $\sigma^\lor\cap M$.
\end{itemize}

There are monoid isomorphisms $\sigma\cap N\cong\mathrm{Hom}(\sigma^\lor\cap M,\mathbb{N})$ and
\[\mathrm{Hom}(\sigma^\lor\cap M,\mathbb{N})\to \operatorname{Ker}A\cap\mathbb{N}^n~;~\varphi\mapsto \begin{pmatrix} \varphi(u_1) \\ \vdots \\ \varphi(u_n)\end{pmatrix}.\]
Combining them we have a monoid isomorphism $\sigma\cap N\cong \operatorname{Ker}A\cap\mathbb{N}^n$ given by
\[v\mapsto \langle\cdot,v\rangle\mapsto (\langle u_j,v\rangle)_{j=1,2,\dots,n}.\]
Since $\{v_1,\dots,v_m\}$ is a Hilbert basis of $\sigma\cap N$, the column vectors $(\langle u_j,v_i\rangle)_{j=1,2,\dots,n}~(i=1,2,\dots,m)$ form the Hilbert basis of 
$\operatorname{Ker}A\cap\mathbb{N}^n$. Hence we have $B=(\langle u_j,v_i\rangle)_{j,i}$.

Now let $\tau\preceq\sigma$ be a face and $Y'\subset Y$ be the corresponding normal subvariety. 
We denote by $B'$ the twisting matrix for the (quotient) cone $q(\sigma)\subset N(\tau)_\mathbb{R}$.

\begin{prop}\label{prop:PBQ}
There exist matrices $P\in M_{n',n}(\mathbb{N}),~Q\in M_{m,m'}(\mathbb{N})$ such that $B'=PBQ$.
\end{prop}
\begin{proof}
Let $\overline{v_i}:=q(v_i)\in N(\tau)$. Since $\overline{v_1},\dots,\overline{v_m}$ generate $q(\sigma)\cap N(\tau)=q(\sigma\cap N)$ we can write
\[(v_1',\dots,v_{m'}')=(\overline{v_1},\dots,\overline{v_m})Q,
\quad Q\in M_{m,m'}(\mathbb{N}).\]
Noting that $u_1,\dots,u_{n'}$ annihilate $\tau$ we see that
\begin{align*}
B' &= \begin{pmatrix}
u_1 \\ \vdots \\ u_{n'}
\end{pmatrix}(v_1',\dots,v_{m'}') 
= \begin{pmatrix}
u_1 \\ \vdots \\ u_{n'} \end{pmatrix}
(\overline{v_1},\dots,\overline{v_m})Q\\
&= \begin{pmatrix}
u_1 \\ \vdots \\ u_{n'} \end{pmatrix}
(v_1,\dots,v_m)Q 
= ({\rm Id}_{n'}~O)\begin{pmatrix}
u_1 \\ \vdots \\ u_n \end{pmatrix}
(v_1,\dots,v_m)Q \\
&= PBQ, \qquad P=({\rm Id}_{n'}~O).
\end{align*}
\end{proof}

\begin{thm}\label{thm:tfact}
Let $W$ be a factorial complex space, and let $Y=Y_\sigma\subset\mathbb C^n$ be a normal
affine toric variety containing the origin, in the Hilbert basis embedding of
\S~\ref{subsec:affine}. Then every holomorphic map $f\colon W\to Y$ has a twisted factorization with respect to the twisting matrix $B=(\langle u_j,v_i\rangle)$ and any fixed matrix $E$ whose columns form a $\mathbb{Z}$-basis of $\operatorname{Ker}B$.
\end{thm}
\begin{proof}
First assume that $W$ is connected hence irreducible (see \cite[pp.~124, 168]{GR84}). By \cite[Lemma 2.9]{LL16} and Theorem~\ref{thm:orbit-cone} there is a face $\tau\preceq\sigma$ 
and corresponding subvariety $Y'=V(\tau)$ such that $f(W)\subset Y'\subset\mathbb{C}^{n'}$ and $f\colon W\to Y'$ is nondegenerate.
We adopt the notation introduced at the beginning of \S~\ref{subsec:adv} for the face $\tau$. First we assume that $Y'$ is not $\mathbb{C}^d$ for any $d\geqslant 0$. 
Then by Theorem~\ref{thm:TFT} there are matrices $B'\in M_{n',m'}(\mathbb{N}),~E'\in M_{m',l'}(\mathbb{Z})$ and holomorphic line bundles $N_\lambda'~(\lambda=1,\dots,l')$ 
such that if we set $L'=(N')^{E'}$, there are global sections $\sigma_i'\in H^0(W,L_i')$ such that $f=(\sigma')^{B'}$. 
We then modify $N'_\lambda$'s and $\sigma'_i$'s so that they fit ambient data. Let us define $I=\{i\mid v_i\in \tau\}$ and 
\[\mathbb{Z}^I=\{(a_1,\dots,a_m)\in\mathbb{Z}^m\mid a_i=0,~\forall i\notin I\}.\]
We first explain the lattice maps used in the construction. Define
\[
\iota\colon N\to\mathbb Z^n,\qquad
v\mapsto(\langle u_j,v\rangle)_{j=1}^n.
\]
The elements $u_j$ generate $M$ as a group, so $\iota$ is injective.
The lattice points in a rational cone generate, as a group, the lattice in its real span:
if $w$ is a lattice point in the relative interior, then any lattice point $v$ in the span
is the difference of $v+kw$ and $kw$, both in the cone for sufficiently large $k$.
Consequently, $\operatorname{Im}B=\iota(N)$.
Since $\tau$ is a face, its Hilbert basis is $\{v_i\mid i\in I\}$, and these elements
generate $N_\tau$ as a group. Similarly, $u_1,\dots,u_{n'}$ generate
$M\cap\tau^\perp$ as a group. Thus, identifying $\operatorname{Im}B'$ with $N(\tau)$
by the analogous pairing map, the restriction of $P$ to $\operatorname{Im}B$
corresponds to the quotient map $N\to N/N_\tau$. In particular,
\[
\operatorname{Ker}P\cap\operatorname{Im}B=\iota(N_\tau),\qquad
B(\mathbb Z^I)=\iota(N_\tau),\qquad
P(\operatorname{Im}B)=\operatorname{Im}B'.
\]
The following diagram records these maps; the two lower horizontal rows are exact:
\[
\xymatrix{
& & \mathbb{Z}^I \ar@{->>}[r]^-{B|_{\mathbb{Z}^I}} \rule[-3pt]{0pt}{0pt} \ar@{^{(}->}[d] & \operatorname{Ker}P\cap \operatorname{Im}B \rule[-3pt]{0pt}{0pt}\ar@{^{(}->}[d] & \\
0 \ar[r] & \mathbb{Z}^l \ar[r]^-E & \mathbb{Z}^m \ar[r]^-B & \operatorname{Im}B \ar@{->>}[d]^-{P} \ar[r] & 0 \\
0 \ar[r] & \mathbb{Z}^{l'} \ar[r]^-{E'} & \mathbb{Z}^{m'} \ar[u]_-{Q} \ar[r]^-{B'} & \operatorname{Im}B' \ar[r] & 0
}
\]
Here $P,Q$ are the matrices from Proposition~\ref{prop:PBQ}.
The preceding identifications prove that the displayed restrictions of $B$ and $P$
are well-defined and surjective. Since $N_\tau$ is a primitive sublattice of $N$,
$\operatorname{Ker}P\cap\operatorname{Im}B$ is a direct summand of $\operatorname{Im}B$.
Choose a right inverse of $B|_{\mathbb Z^I}$ with values in $\mathbb Z^I$.
On a complementary summand of $\operatorname{Im}B$, choose lifts under $B$ of a basis;
this extends the chosen right inverse to a homomorphism
$\widetilde B\colon\operatorname{Im}B\to\mathbb Z^m$ satisfying $B\widetilde B=\mathrm{Id}$.
For $a\in\operatorname{Ker}B'$, the relation $B'=PBQ$ gives
$BQa\in\operatorname{Ker}P\cap\operatorname{Im}B$. Hence
$(Q-\widetilde B BQ)a\in\operatorname{Ker}B$, and there is a unique matrix
$R\in M_{l,l'}(\mathbb Z)$ such that
\[
ER=(Q-\widetilde B BQ)E'.
\]

We define $N=(N')^R$ and $L=N^E=(N')^{ER}$. Since $\operatorname{Im}(BQE')\subset\operatorname{Ker}P\cap\operatorname{Im}B$,
the columns of $\widetilde B BQE'$ belong to $\mathbb Z^I$.
Thus, for $i\notin I$, the $i$th row of $ER$ equals the corresponding row of $QE'$. Therefore we can define a holomorphic section of $L_i$ by setting
\[
\sigma_i=\begin{cases}
(\sigma_1')^{q_{i1}}\otimes\dots\otimes(\sigma'_{m'})^{q_{im'}}, & i\notin I,\\
0, & i\in I.
\end{cases}\eqno{(\ast)}
\]
We now check that the equality $f=\sigma^B$ holds. For $j=1,2,\dots,n'$ and $i\in I$ we have $b_{ji}=0$ and therefore
\[
\sigma_1^{b_{j1}}\otimes\dots\otimes\sigma_m^{b_{jm}}=\bigotimes_{i\notin I}\sigma_i^{b_{ji}}
=\bigotimes_{k=1}^{m'}(\sigma_k')^{\sum_{i\notin I} b_{ji}q_{ik}}=\bigotimes_{k=1}^{m'}(\sigma_k')^{b_{jk}'}=f_j.
\]
Here we use the convention $0^0=1$.
For $j>n'$, the element $u_j$ does not vanish identically on $\tau$.
Since the $v_i$ with $i\in I$ generate $\tau\cap N$, there is an $i\in I$
such that $b_{ji}=\langle u_j,v_i\rangle>0$.
Thus $\bigotimes_{i=1}^m\sigma_i^{b_{ji}}=0=f_j$.

If $Y'=\mathbb{C}^d$ for some $d\geqslant0$, we can take $B'=I_d$,
so $l'=0$. We set $L_i'=\mathbb{I}_W$ and
$\sigma_i'=f_i~(i=1,2,\dots,d)$, with the convention that an empty
tensor product is the trivial line bundle.
Clearly $f=(\sigma')^{B'}$. The same construction in $(\ast)$ also applies in this case.

In general we argue componentwise. Let $W=\bigsqcup_\alpha W_\alpha$ where $W_\alpha$ is a connected component of $W$. 
By the above argument there are line bundles $N_{\lambda,\alpha}~(\lambda=1,\dots,l)$ and sections $\sigma_{i,\alpha}~(i=1,\dots,m)$ such that $f|_{W_\alpha}=\sigma_\alpha^B$. 
Now we define line bundles on $W$ by setting $N_\lambda|_{W_\alpha}=N_{\lambda,\alpha}$ and sections $\sigma_i$ on $W$ by setting $\sigma_{i}|_{W_\alpha}=\sigma_{i,\alpha}$. 
Then we have $f=\sigma^B$.
\end{proof}

We give two corollaries. The first one is an extension theorem.

\begin{cor}\label{cor:ext}
Let $Y$ be a normal affine toric variety. Let $X'$ be a closed factorial subvariety of a reduced Stein space $X$ such that $H^q(X,\mathbb{Z})\to H^q(X',\mathbb{Z})$
is surjective for $q=1,2$. Then every holomorphic map $X'\to Y$ extends to a holomorphic map $X\to Y$.
\end{cor}
\begin{proof}
By the decomposition obtained in the proof of Proposition~\ref{prop:0inY},
there is a normal affine toric variety $Y'\subset\mathbb{C}^{n'}$ 
containing $0\in\mathbb{C}^{n'}$ such that $Y=Y'\times(\mathbb{C}^\ast)^s$. 
We may write $f=(f_1,f_2)$ according to the decomposition. By Theorem~\ref{thm:tfact} $f_1\colon X'\to Y'$ has a twisted factorization. 
By the assumption that $H^2(X,\mathbb{Z})\to H^2(X',\mathbb{Z})$ is surjective, all line bundles on $X'$ extend to $X$. 
We first extend each $N_\lambda$ on $X'$ to $\widetilde{N}_\lambda$ on $X$ and put $\widetilde{L}=\widetilde{N}^E$.
By Cartan's Theorem~B (\cite[Theorem 7.4.3]{Hor90}, \cite[Theorem 4.4.5]{GR79}) there is an extension $\widetilde\sigma_i\in H^0(X,\widetilde{L}_i)$ of each $\sigma_i$. 
Define $F_1=\widetilde\sigma^B$. This is a map $X\to\mathbb{C}^{n'}$ since $BE=0$ and moreover a map $X\to Y'$ since $AB=0$.
By the surjectivity of $H^1(X,\mathbb{Z})\to H^1(X',\mathbb{Z})$ and elementary obstruction theory, $f_2\colon X'\to (\mathbb{C}^\ast)^s$ extends continuously to all of $X$. 
Since $(\mathbb{C}^\ast)^s$ is an Oka manifold, there is a holomorphic extension $F_2$ of $f_2$. Finally we put $F=(F_1,F_2)\colon X\to Y'\times(\mathbb{C}^\ast)^s$.
\end{proof}

The second one is an approximation theorem.

\begin{cor}
Let $Y$ be a normal affine toric variety. Let $\Omega$ be a Stein Runge open subset of a factorial Stein space $X$ such that $H^q(X,\mathbb{Z})\to H^q(\Omega,\mathbb{Z})$ 
is surjective for $q=1,2$. Then every holomorphic map $\Omega\to Y$ can be approximated uniformly on compact subsets of $\Omega$ by holomorphic maps $X\to Y$.
\end{cor}
\begin{proof}
Let $f\colon\Omega\to Y$ be holomorphic. Write $Y=Y'\times(\mathbb C^*)^s$, and write $f=(f_1,f_2)$ as before.
By Theorem~\ref{thm:tfact}, the map $f_1$ has a twisted factorization $f_1=\sigma^B$,
with $L=N^E$. 
The surjectivity of $H^2(X,\mathbb Z)\to H^2(\Omega,\mathbb Z)$ allows us
to extend each $N_\lambda$ to a holomorphic line bundle $\widetilde N_\lambda$ on $X$.
Put $\widetilde L=\widetilde N^E$ and identify $\widetilde L_i|_\Omega$ with $L_i$.
Runge approximation for coherent analytic sheaves gives sections
$\sigma_{i,j}\in H^0(X,\widetilde L_i)$ converging to $\sigma_i$ uniformly on compact
subsets of $\Omega$, with respect to fixed hermitian metrics.
Then $F_{1,j}=\sigma_j^B$ are holomorphic maps $X\to Y'$ converging to $f_1$
uniformly on compact subsets of $\Omega$.

For the torus component, let $\delta_X$ and $\delta_\Omega$ denote the connecting
homomorphisms of the exponential sequences on $X$ and $\Omega$, respectively.
Since $X$ is Stein, Cartan's Theorem B shows that $\delta_X\colon\mathcal O^*(X)\to H^1(X,\mathbb Z)$ is surjective. 
By the surjectivity of $H^1(X,\mathbb Z)\to H^1(\Omega,\mathbb Z)$ and naturality
of the connecting homomorphisms, we can therefore choose $g\in\mathcal O(X,(\mathbb C^*)^s)$ such that
$\delta_\Omega(g|_\Omega)=\delta_\Omega(f_2)$, componentwise.
Exactness of the exponential sequence on $\Omega$ then gives
$h\in\mathcal O(\Omega)^s$ satisfying $f_2/g|_\Omega=e^h$.
Since $\Omega$ is Runge in $X$, choose $h_j\in\mathcal O(X)^s$ converging to $h$
uniformly on compact subsets of $\Omega$.
The holomorphic maps $F_{2,j}=g e^{h_j}\colon X\to(\mathbb C^*)^s$ then converge to $f_2$
uniformly on compact subsets of $\Omega$.
Finally, $F_j=(F_{1,j},F_{2,j})\colon X\to Y$ gives the required approximation.
\end{proof}


\section{Simultaneous approximation and interpolation theorem}\label{sec:BOPAI}

In this section we prove a simultaneous approximation and interpolation theorem. 
The main idea of the proof is already contained in \cite[Theorem 3.4.1]{For17}. 
To prove this we have to compare two twisted factorizations of a holomorphic map. 
The following lemma gives a condition under which two twisted factorizations are isomorphic.

\begin{lem}\label{lem:comp}
Let $X$ be a complex manifold, and let $Y=Y_\sigma\subset\mathbb C^n$ be a normal affine
toric variety containing the origin, in the Hilbert basis embedding of \S~\ref{subsec:affine}.
Let $f\colon X\to Y$ be a nondegenerate holomorphic map. Assume that 
\[\operatorname{codim}_Xf^{-1}(Y_{\rm sing})\geqslant 2.\] 
Then for twisted factorization data $(N_\lambda,\sigma_i)$ and $(N_\lambda',\sigma_i')$ 
of $f$ with respect to the matrix $B=(\langle u_j,v_i\rangle)$ and any fixed matrix $E$, 
there is a (holomorphic) isomorphism $N_\lambda\cong N_\lambda'$ 
such that under the induced isomorphism $L_i\cong L_i'$ we have $\sigma_i=\sigma_i'$.
\end{lem}
\begin{proof}
It is enough to show that for all $i=1,2,\dots,m$ we have $\div(\sigma_i)=\div(\sigma_i')$. Indeed, recall the exact sequence (cf.~\cite[Lemma 2.1]{LL16}):
\[\xymatrix{
0 \ar[r] & \mathbb{Z}^l \ar[r]^-E & \mathbb{Z}^m \ar[r]^-B & \mathbb{Z}^n \ar[r]^-C & \mathbb{Z}^p \ar[r] & 0.
}\]
Since $\operatorname{Im}B$ is a direct summand in $\mathbb{Z}^n$,
tensoring with $\mathcal{O}^\ast(X)$ preserves exactness:
\[\xymatrix{
0 \ar[r] & \mathcal{O}^\ast(X)^l \ar[r]^-E & \mathcal{O}^\ast(X)^m \ar[r]^-B & \mathcal{O}^\ast(X)^n \ar[r]^-C & \mathcal{O}^\ast(X)^p \ar[r] & 0.
}\eqno{(\ast\ast)}\]
If we have $\div(\sigma_i)=\div(\sigma_i')$ for all $i$, 
then $\psi_i=\sigma_i'/\sigma_i$ gives a nonvanishing section 
of the line bundle $L_i^{-1}\otimes L_i'$, i.e., $\psi_i$ gives
an isomorphism $L_i\cong L_i'$. 
Choose a left inverse $\widetilde{E}=(\widetilde{e}_{\lambda i})\in M_{l,m}(\mathbb{Z})$ of $E$ and put
\[\varphi_\lambda=\psi_1^{\widetilde{e}_{\lambda 1}}\otimes\dots\otimes\psi_m^{\widetilde{e}_{\lambda m}}\in H^0(X,N_\lambda^{-1}\otimes N_\lambda').\]
Then $\varphi_\lambda\colon N_\lambda \to N_\lambda'$ is also an isomorphism. 
Let us consider another isomorphism $(\varphi^E)_i\colon L_i\to L_i'$. 
The quotient $r_i:=(\varphi^E)_i/\psi_i$ is a nonvanishing holomorphic function on $X$. 
We can easily see that $r^B=1$. 
By the exactness of ($\ast\ast$) there is $\eta\in\mathcal{O}^\ast(X)^l$ such that $r=\eta^E$. 
Here we have $\eta=\eta^{\widetilde{E}E}=(\eta^E)^{\widetilde{E}}=r^{\widetilde{E}}$.
Again we can see that $r^{\widetilde{E}}=1$ so $\eta=1$ and therefore $r=1$, i.e., $\varphi^E=\psi$.

Now we prove that $\div(\sigma_i)=\div(\sigma_i')$ for all $i$. Let us write
\[\div(\sigma_i)=\sum_{D\,:\,\text{prime divisor}}a_{iD}D.\]
We want to show that the coefficients $a_{iD}$ are uniquely determined by $f$. 
For a prime divisor $D$ we associate a lattice point $\nu_D\in\sigma\cap N$ as follows: The mapping
\[u\in\sigma^\lor\cap M\mapsto \operatorname{ord}_D(f^\ast\chi^u)\in\mathbb{N}\]
is a morphism of monoids. 
By the fundamental monoid isomorphism $\mathrm{Hom}(\sigma^\lor\cap M,\mathbb{N})\cong\sigma\cap N$, 
there is a unique corresponding lattice point $\nu_D\in\sigma\cap N$. In other words $\nu_D$ is the unique point satisfying
\[\langle u,\nu_D\rangle=\operatorname{ord}_D(f^\ast\chi^u),\quad\forall u\in\sigma^\lor\cap M.\]
For the Hilbert basis $\{u_j\}$ of $\sigma^\lor\cap M$ we have
\[\langle u_j,\nu_D\rangle=\operatorname{ord}_D(f_j)=\sum_ib_{ji}a_{iD}=\left\langle u_j,\sum_i a_{iD}v_i\right\rangle.\]
Since the $u_j$ generate $\sigma^\lor\cap M$, we have $\nu_D=\sum_i a_{iD}v_i$.
Take the maximal face $\tau_D\preceq\sigma$ such that $f(D)\subset V(\tau_D)$;
then $f(D)\cap O(\tau_D)\ne\emptyset$.
The codimension assumption and Proposition~\ref{prop:smoothsing} imply that $\tau_D$ is smooth.
For every $u\in\sigma^\lor\cap\tau_D^\perp\cap M$, the restriction
$(f^*\chi^u)|_D$ is not identically zero, and hence
\[
\langle u,\nu_D\rangle=\operatorname{ord}_D(f^*\chi^u)=0.
\]
Since $\nu_D\in\sigma$, the duality of faces gives $\nu_D\in\tau_D$.

Now $\nu_D=\sum_i a_{iD}v_i$ with $a_{iD}\ge0$, and $\tau_D$ is a face of $\sigma$.
Therefore every $v_i$ occurring with $a_{iD}>0$ belongs to $\tau_D$.
The intersection of the Hilbert basis of $\sigma\cap N$ with $\tau_D$
is the Hilbert basis of $\tau_D\cap N$.
Because $\tau_D$ is smooth, this basis consists of its primitive ray generators,
which are linearly independent. Thus the coefficients $a_{iD}$ are uniquely
determined by $\nu_D$, and hence by $f$.
Applying the same argument to $\sigma_i'$ proves
$\div(\sigma_i)=\div(\sigma_i')$ for every $i$.
\end{proof}

\begin{thm}\label{thm:simAI}
Let $X$ be a Stein manifold, $X'\subset X$ be a closed submanifold, 
and $K\subset X$ be an $\mathcal{O}(X)$-convex compact set. 
Let $Y$ be a normal affine toric variety.
Let $f\colon X\to Y$ be a continuous map which is holomorphic 
on $X'$ and in some neighborhood $U$ of $K$. Assume that
\begin{enumerate}[label={\rm (\arabic*)}]
\item $f$ is nondegenerate on $X'\cap U$ and $\operatorname{codim}_{X'\cap U}f^{-1}(Y_{\rm sing})\geqslant2;$
\item the restriction map $H^2(X,\mathbb{Z})\to \check{H}^2(X'\cup K,\mathbb{Z})$ is surjective.
\end{enumerate}
Then there is a sequence of holomorphic maps $F_j\colon X\to Y,~j=1,2,\dots$ such that $F_j|_{X'}=f|_{X'}$ and $F_j\to f$ uniformly on $K$ as $j\to\infty$.
Moreover, for each $j$ there is a homotopy $f_{j,\theta}\colon X\to Y$ $(0\leqslant\theta\leqslant1)$ 
with $f_{j,0}=f$ and $f_{j,1}=F_j$ such that, for every $\theta$, 
$f_{j,\theta}|_{X'}=f|_{X'}$ and $f_{j,\theta}$ is holomorphic on a neighborhood of $K$, and
$f_{j,\theta}\to f$ as $j\to\infty$ uniformly on $K$ and uniformly in $\theta\in[0,1]$.
\end{thm}

\begin{proof}
By the decomposition obtained in the proof of Proposition~\ref{prop:0inY}, 
there exists a normal affine toric variety $Y'\subset\mathbb{C}^{n'}$ containing $0\in\mathbb{C}^{n'}$ 
such that $Y=Y'\times(\mathbb{C}^\ast)^s$. 
We may write $f=(f_1,f_2)$ according to the decomposition. 
Since $f_2\colon X\to(\mathbb{C}^\ast)^s$ is continuous, 
$f_2|_{X'}$ and $f_2|_U$ are holomorphic, 
and $(\mathbb{C}^\ast)^s$ is an Oka manifold, 
the Oka principle (\cite[Theorem 5.4.4]{For17}) gives a sequence of holomorphic maps 
$F_{2,j}\colon X\to(\mathbb{C}^\ast)^s,~j=1,2,\dots$, together with homotopies 
$f_{2,j,\theta}\colon X\to(\mathbb{C}^\ast)^s$ $(0\leqslant\theta\leqslant1)$ 
from $f_{2,j,0}=f_2$ to $f_{2,j,1}=F_{2,j}$, 
such that $f_{2,j,\theta}|_{X'}=f_2|_{X'}$, each $f_{2,j,\theta}$ is holomorphic
on a neighborhood of $K$, and $f_{2,j,\theta}\to f_2$ as $j\to\infty$ uniformly 
on $K$ and uniformly in $\theta$. By Theorem~\ref{thm:tfact} $f_1|_{X'}$ and 
$f_1|_U$ have twisted factorizations $(N_{\lambda,X'},\sigma_{i,X'})$ and 
$(N_{\lambda,U},\sigma_{i,U})$. 
Since
\[(\sigma_{X'}|_{X'\cap U})^B=f_1|_{X'\cap U}=(\sigma_{U}|_{X'\cap U})^B\]
and $\operatorname{codim}_{X'\cap U}f^{-1}(Y_{\rm sing})\geqslant2$, 
Lemma~\ref{lem:comp} ensures that there exist isomorphisms 
\[\varphi_\lambda\colon N_{\lambda,U}|_{X'\cap U}\to N_{\lambda,X'}|_{X'\cap U}\] 
such that the induced isomorphisms
\[
(\varphi^E)_i\colon L_{i,U}|_{X'\cap U}\to L_{i,X'}|_{X'\cap U}
\]
satisfy
\[
(\varphi^E)_i(\sigma_{i,U}|_{X'\cap U})=\sigma_{i,X'}|_{X'\cap U}.
\]
Now choose a Stein neighborhood $U_1$ of $K$ such that $K\subset U_1\Subset U$ and 
put $Z=X'\cup\overline{U_1}$ (see figure~\ref{fig:sets} below). 
Let us consider its open covering $\{V_1=Z\cap U,~V_2=X'\setminus \overline{U_1}\}$.

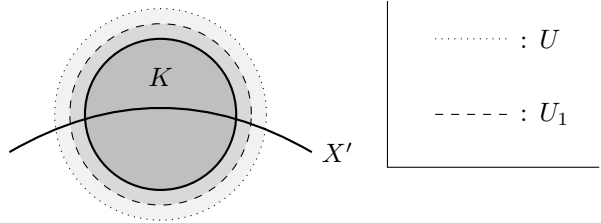
\begin{figure}[h]
\begin{tikzpicture}
\filldraw[dotted,fill=black!5] (0,0) circle (1.4);
\filldraw[dashed,fill=black!15] (0,0) circle (1.2);
\filldraw[thick,draw=black,fill=black!25] (0,0) circle (1);
\node at (0,.5) {$K$};
\draw[thick] (-2,-.5) to[bend left] (2,-.5) node[right]{$X'$};

\draw (3,1.5)--(3,-.7)--(5.8,-.7);
\node[right] at (3.5,1) {\tikz{\draw[dotted](0,0)--(1,0) node[right]{: $U$};}};
\node[right] at (3.5,0) {\tikz{\draw[dashed](0,0)--(1,0) node[right]{: $U_1$};}};
\end{tikzpicture}
\caption{$X'$, $K$, $U$, and $U_1$}
\label{fig:sets}
\end{figure}

We then get topological vector bundles $N_{\lambda,Z}$ by gluing
$N_{\lambda,X'}$ and $N_{\lambda,U}$ via $\varphi_\lambda$ 
on $V_1\cap V_2=X'\cap(U\setminus\overline{U_1})$. 
We denote by $c_1(N_{\lambda,Z})\in H^2(Z,\mathbb{Z})$ its first Chern class. 
Choose relatively compact open analytic polyhedra $U_k$ for $k\ge2$
such that $K\subset U_k\Subset U_1$, $\overline{U_{k+1}}\subset U_k$, and
$K=\bigcap_{k\ge2}\overline{U_k}$. Put $Z_k=X'\cup\overline{U_k}$.
By continuity of \v{C}ech cohomology for a decreasing sequence
of compact Hausdorff spaces
(\cite[Theorems 6.6.5, 6.6.6, and Corollary 6.8.8]{Spa94}), we have
\[
\check{H}^q(K,\mathbb{Z})
=\varinjlim H^q(\overline{U_k},\mathbb{Z}),\qquad
\check{H}^q(X'\cap K,\mathbb{Z})
=\varinjlim H^q(X'\cap\overline{U_k},\mathbb{Z}).
\]
The Mayer--Vietoris exact sequence and the Five Lemma then give
\[
\check{H}^2(X'\cup K,\mathbb{Z})
=\varinjlim H^2(Z_k,\mathbb{Z}).
\]
By the assumption (2) of this theorem, there is an element
$x_\lambda\in H^2(X,\mathbb{Z})$ such that 
$x_\lambda|_{X'\cup K}=c_1(N_{\lambda,Z})|_{X'\cup K}$. 
By the above continuity there is an index $k$ such that $x_\lambda|_{Z_k}=c_1(N_{\lambda,Z})|_{Z_k}$. 
Since there are only finitely many $\lambda$, we may choose the same index $k$ for all of them.
We now replace $U$ by $U_k$ and $Z$ by $Z_k$, restricting $N_{\lambda,Z}$ accordingly.
By the classical Oka principle there is a holomorphic line bundle
$N_{\lambda,X}$ on $X$ such that $c_1(N_{\lambda,X})=x_\lambda$.
Since $c_1(N_{\lambda,X})|_Z=c_1(N_{\lambda,Z})$, there are topological
isomorphisms
\[\alpha_\lambda^{\mathcal C}\colon N_{\lambda,X}|_{X'}\to N_{\lambda,X'},\quad
\beta_\lambda^{\mathcal C}\colon N_{\lambda,X}|_U\to N_{\lambda,U}\]
such that $\varphi_\lambda\circ\beta_\lambda^{\mathcal C}=\alpha_\lambda^{\mathcal C}$
on $X'\cap U$. Again by the Oka principle, since $X'$ and $U$ are Stein,
there are holomorphic isomorphisms
\[\alpha_\lambda^{\mathcal O}\colon N_{\lambda,X}|_{X'}\to N_{\lambda,X'},
\quad\beta_\lambda^{\mathcal O}\colon N_{\lambda,X}|_U\to N_{\lambda,U},\]
which need not be compatible. We modify them so that
$\varphi_\lambda\circ\beta_\lambda^{\mathcal O}=\alpha_\lambda^{\mathcal O}$
holds on $X'\cap U$. The function
$\alpha_\lambda^{\mathcal O}\circ(\alpha_\lambda^{\mathcal C})^{-1}$ is a
continuous map $X'\to\mathbb C^*$, and on a Stein space every such map is
homotopic to a holomorphic one. Hence, after multiplying
$\alpha_\lambda^{\mathcal O}$ by a suitable nowhere vanishing holomorphic
function, and similarly for $\beta_\lambda^{\mathcal O}$, we may assume that
there are continuous functions $a_\lambda\colon X'\to\mathbb C$ and
$b_\lambda\colon U\to\mathbb C$ such that
\[\alpha_\lambda^{\mathcal O}=e^{a_\lambda}\alpha_\lambda^{\mathcal C},\quad
\beta_\lambda^{\mathcal O}=e^{b_\lambda}\beta_\lambda^{\mathcal C}.\]
Then on $X'\cap U$ we have
\[(\varphi_\lambda\circ\beta_\lambda^{\mathcal O})\circ(\alpha_\lambda^{\mathcal O})^{-1}
=e^{b_\lambda-a_\lambda}(\varphi_\lambda\circ\beta_\lambda^{\mathcal C})\circ(\alpha_\lambda^{\mathcal C})^{-1}
=e^{b_\lambda-a_\lambda}.\]
It follows that
$(b_\lambda-a_\lambda)|_{X'\cap U}\in\mathcal O(X'\cap U)$. Since $U$ is Stein
and $X'\cap U$ is a closed submanifold of $U$, Cartan's Theorem B gives a
holomorphic function $\eta_\lambda\in\mathcal O(U)$ extending
$(b_\lambda-a_\lambda)|_{X'\cap U}$. Finally, we set
$\alpha_\lambda:=\alpha_\lambda^{\mathcal O}$ and
$\beta_\lambda:=e^{-\eta_\lambda}\beta_\lambda^{\mathcal O}$. These are
holomorphic isomorphisms satisfying
$\varphi_\lambda\circ\beta_\lambda=\alpha_\lambda$ on $X'\cap U$.
Let us put $L_X=N_X^E$ and define
\[s_{i,X'}=(\alpha^E)^{-1}_i(\sigma_{i,X'})\in H^0(X',L_{i,X}|_{X'}),
\quad s_{i,U}=(\beta^E)_i^{-1}(\sigma_{i,U})\in H^0(U,L_{i,X}).\]
Since $\varphi_\lambda\circ\beta_\lambda=\alpha_\lambda$ and 
$(\varphi^E)(\sigma_U)=\sigma_{X'}$ on $X'\cap U$ we see that $s_{i,X'}|_{X'\cap U}=s_{i,U}|_{X'\cap U}$. 
By Cartan's Theorem B there is a holomorphic section 
$s_{i,X}\in H^0(X,L_{i,X})$ such that $s_{i,X}|_{X'}=s_{i,X'}$. Note that
\[s_{i,U}-s_{i,X}|_U\in H^0(U,\mathcal{I}_{X'}\otimes\mathcal{O}_X(L_{i,X})).\]
By the Oka--Weil approximation for coherent analytic sheaves (\cite[Theorem 2.6.8]{For17}) 
there is a sequence of global holomorphic sections 
$t_{ij,X}\in H^0(X,\mathcal{I}_{X'}\otimes\mathcal{O}_X(L_{i,X}))$ such that, 
for some hermitian metric on $L_{i,X}$, 
we have $t_{ij,X}\to s_{i,U}-s_{i,X}~(j\to\infty)$ uniformly on $K$. 
Since $t_{ij,X}$ vanishes on $X'$ we see that $(s_{i,X}+t_{ij,X})|_{X'}=s_{i,X'}$. 
So if we put
\[F_{1,j}=(s_{X}+t_{j,X})^B,\]
then $F_{1,j}|_{X'}=f_1|_{X'}$ and $F_{1,j}\to f_1$ uniformly on $K$.
The sequence of global holomorphic maps $F_j=(F_{1,j},F_{2,j})$
satisfies the first part of the conclusion.

It remains to construct the homotopies. Put $U'=U_{k+1}$, so that $K\subset U'$ and
$\overline{U'}\subset U$, and choose a continuous function $\kappa\colon X\to[0,1]$
with compact support in $U$ which equals $1$ on a neighborhood of $\overline{U'}$.
The continuous sections
\[s_i^c=s_{i,X}+\kappa\,(s_{i,U}-s_{i,X})\]
of $L_{i,X}$ satisfy $s_i^c|_{X'}=s_{i,X'}$ and
$s_i^c=s_{i,U}$ near $\overline{U'}$. Since $BE=0$ and $AB=0$, $(s^c)^B$ is a continuous
map $X\to Y'$, and it agrees with $f_1$ on $X'\cup\overline{U'}$.
As $Y'$ is contractible and $(X,X'\cup\overline{U'})$ is a CW pair (cf.~Remark~\ref{rem:Cech}),
$f_1$ and $(s^c)^B$ are homotopic relative to $X'\cup\overline{U'}$. Next,
\[\theta\mapsto\bigl(s_X+(1-\theta)\,\kappa\,(s_U-s_X)+\theta\,t_{j,X}\bigr)^B
\quad(0\leqslant\theta\leqslant1)\]
joins $(s^c)^B$ to $F_{1,j}$. This homotopy is fixed on $X'$ and holomorphic on $U'$
for every $\theta$, since $\kappa=1$ there. On $K$ the sections inside the parentheses differ
from $s_U$ by $\theta\bigl(t_{j,X}-(s_U-s_X)\bigr)$, so the homotopy converges to $f_1$
uniformly on $K$ and uniformly in $\theta$ as $j\to\infty$.
Concatenating these two homotopies and combining the result with $f_{2,j,\theta}$,
we obtain $f_{j,\theta}$. 
\end{proof}

\begin{rem}\label{rem:Cech}
We use the \v{C}ech cohomology $\check{H}^\bullet$ because it has the continuity property. 
Note that $X'$, $\overline{U_k}$, and $Z_k=X'\cup\overline{U_k}$ are closed semi-analytic subsets of $X$, 
so $X$ admits a triangulation in which they are subcomplexes (\cite[Theorem 2, p.~464]{Loj64}); 
in particular, $(X,Z_k)$ is a CW pair, and singular and \v{C}ech cohomology agree on these spaces. 
But the space $X'\cup K$ may not be even locally contractible 
(recall that any compact set $K\subset\mathbb{R}^n\subset\mathbb{C}^n$ is 
$\mathcal{O}(\mathbb{C}^n)$-convex).
\end{rem}

We conclude with a counterexample showing that the codimension hypothesis in
Theorem~\ref{thm:simAI} cannot be omitted, even when $X'\cup K$ is contractible.
The prescribed data are holomorphic on $X'$ and near $K$, but admit no common holomorphic
extension to a neighborhood of $X'\cup K$. This distinction is essential: Gromov's definition
of the $\mathrm{Ell}_\infty$-property requires the initial map to be holomorphic on such a
neighborhood (\cite[\S 3.1.A, p.~878]{Gro89}). Thus the following example demonstrates an
obstruction to simultaneous approximation and interpolation in the setting considered here;
it does not by itself disprove the $\mathrm{Ell}_\infty$-property.

\begin{ex}\label{ex:notAI}
Our setting is as follows:
\begin{itemize}
\item $X=\mathbb{C}^3_{xyz}$, $X'=\{z=0\}\cong\mathbb{C}^2_{xy}$, 
$K=\{|x|\leqslant 2,~|y-1|\leqslant 1/2,~|z|\leqslant 1\}$.
\item $Y=\{w^2=uv\}$.
\item There is a single-valued holomorphic branch of $\sqrt{y}$ 
in a neighborhood of $K$. We choose the branch satisfying $\sqrt{1}=1$. 
Let us put
\begin{alignat*}{3}
p(x,y,z)&=x-\sqrt{y},&q(x,y,z)&=x+\sqrt{y},&h(x,y,z)&=pq=x^2-y, \\
a(x,y,z)&=p(x,y,z)+z,~&b(x,y,z)&=p(x,y,z)-z,~&c(x,y,z)&=q(x,y,z)^2+z.
\end{alignat*}
Consider the continuous map
\[
f=
\begin{cases}
(h^2,h^2,h^2) & \text{on}~X', \\
(a^2c,b^2c,abc) & \text{on~a~neighborhood~of}~K.
\end{cases}
\]
\end{itemize}
The two formulas agree on the overlap: when $z=0$, we have $a=b=p$ and $c=q^2$,
so $(a^2c,b^2c,abc)=(h^2,h^2,h^2)$.
For completeness, we construct a continuous extension explicitly.
Let $D=\{y\in\mathbb C:|y-1|<1\}$, and choose a continuous function
$\chi\colon\mathbb C\to[0,1]$ with compact support in $D$ which equals $1$
on a neighborhood of $\{|y-1|\le1/2\}$.
On $\mathbb C\times D\times\mathbb C$, use the chosen branch of $\sqrt y$ and set
\[
\widehat a=p+\chi(y)z,\qquad
\widehat b=p-\chi(y)z,\qquad
\widehat c=q^2+\chi(y)z.
\]
Define $f=(\widehat a^2\widehat c,\widehat b^2\widehat c,
\widehat a\widehat b\widehat c)$ there and $f=(h^2,h^2,h^2)$ for $y\notin D$,
where $h=x^2-y$ is defined globally.
Since $\chi$ vanishes near the boundary of $D$, these formulas give a continuous
map $X\to Y$. It agrees with the prescribed map on $X'$ and is holomorphic
on a neighborhood of $K$.

Moreover, the homotopy
\[
r_t(x,y,z)=(x,y,(1-t)z),\qquad 0\le t\le1,
\]
is a strong deformation retraction of $X'\cup K$ onto $X'$.
Thus $X'\cup K$ is contractible, so
$\check H^2(X'\cup K,\mathbb Z)=0$ and the cohomological hypothesis of
Theorem~\ref{thm:simAI} is satisfied.
On the overlap, $f|_{X'}$ is nondegenerate and
$f^{-1}(Y_{\rm sing})\cap X'=\{h=0\}=\{y=x^2,z=0\}$,
which has codimension one. Hence this example violates precisely the
codimension hypothesis of that theorem.

Assume that there is a sequence of holomorphic maps $F_j\colon X\to Y$ 
such that $F_j|_{X'}=f|_{X'}$ and $\|f-F_j\|_K\leqslant 1/j$. 
Since $\operatorname{Pic}(X)=\{\mathbb{I}_X\}$, a twisted factorization of $F_j$ has the form
\[F_j=(A_j^2C_j,B_j^2C_j,A_jB_jC_j),\quad A_j,B_j,C_j\in\mathcal{O}(X).\]
Consider the restriction $F_j|_{X'}=f|_{X'}$. Then we have
\[A_j^2C_j=B_j^2C_j=A_jB_jC_j=h^2~\text{on}~X'.\] 
Since $h^2$ is not identically zero on $X'$, the displayed identities imply
$A_j|_{X'}=B_j|_{X'}$. The function $h|_{X'}=x^2-y$ has a reduced irreducible
zero divisor, so
$2\div(A_j|_{X'})+\div(C_j|_{X'})=2\div(h|_{X'})$
gives exactly two possibilities:
\begin{enumerate}[label=(\roman*)]
\item $A_j|_{X'}$ and $B_j|_{X'}$ are nowhere zero and $(C_j|_{X'})=2(h|_{X'})$.
\item $C_j|_{X'}$ is nowhere zero and $(A_j|_{X'})=(B_j|_{X'})=(h|_{X'})$.
\end{enumerate}

\medskip
\textbf{Case (i)}: Since $A_j|_{X'}=B_j|_{X'}$ there is a 
holomorphic function $H_j$ such that $A_j-B_j=zH_j$. We have
\[F_{j,1}+F_{j,2}-2F_{j,3}=C_j(A_j-B_j)^2=z^2C_jH_j^2.\]
Since $C_j(1,1,0)=0$, we obtain
\[
\frac{\del^2}{\del z^2}(F_{j,1}+F_{j,2}-2F_{j,3})(1,1,0)=0.
\]
On the other hand, near $K$ we have
$f_1+f_2-2f_3=c(a-b)^2=4z^2(q^2+z)$.
As $q(1,1,0)=2$, it follows that
\[
\frac{\del^2}{\del z^2}(f_1+f_2-2f_3)(1,1,0)=32.
\]

\medskip
\textbf{Case (ii)}: Here $A_j(-1,1,0)=0$, so
$\del(A_j^2C_j)/\del z$ vanishes at this point.
For the local map, $p(-1,1,0)=-2$ and $q(-1,1,0)=0$, so
$f_1(-1,1,z)=(-2+z)^2z$ near $z=0$. Thus
\[
\frac{\del F_{j,1}}{\del z}(-1,1,0)=0,\qquad
\frac{\del f_1}{\del z}(-1,1,0)=4.
\]

\medskip
Since uniform convergence on $K$ implies convergence of derivatives at interior points of $K$ (Cauchy estimates), and one of the cases (i), (ii) occurs for infinitely many $j$, passing to a subsequence yields a contradiction in either case.
\end{ex}


\section*{Acknowledgements}

The author thanks his supervisor Yuta Kusakabe for his enthusiastic support and insightful comments.
He also thanks Takayuki Koike, his undergraduate supervisor, for his warm encouragement and continued support.


\section*{AI disclosure}

In developing the results presented in this paper, 
the author engaged in iterative discussions with ChatGPT 6 Astra and Claude Opus 5.5
to explore possible approaches, identify gaps in arguments, 
and refine proofs. He critically evaluated the suggestions arising from these
discussions and revised the arguments accordingly. The author has checked the
mathematical content and takes full responsibility for the final version of the paper.



\begin{thebibliography}{CLS11}

\bibitem[CLS11]{CLS11} D.~A. Cox, J.~B. Little and H. Schenck, {\it Toric varieties}, Graduate Studies in Mathematics, 124, Amer. Math. Soc., Providence, RI, 2011; \href{https://mathscinet.ams.org/mathscinet/article?mr=2810322}{MR2810322}

\bibitem[For09]{For09} F. Forstneri\v c, Oka manifolds, C. R. Math. Acad. Sci. Paris {\bf 347} (2009), no.~17-18, 1017--1020; \href{https://mathscinet.ams.org/mathscinet/article?mr=2554568}{MR2554568}

\bibitem[For17]{For17} F. Forstneri\v c, {\it Stein manifolds and holomorphic mappings}, second edition, 
Ergebnisse der Mathematik und ihrer Grenzgebiete. 3. Folge. A Series of Modern Surveys in Mathematics, 56, Springer, Cham, 2017; \href{https://mathscinet.ams.org/mathscinet/article?mr=3700709}{MR3700709}


\bibitem[Gra58]{Gra58} H. Grauert, Analytische Faserungen \"uber holomorph-vollst\"andigen R\"aumen, Math. Ann. {\bf 135} (1958), 263--273; \href{https://mathscinet.ams.org/mathscinet/article?mr=98199}{MR0098199}

\bibitem[GR79]{GR79} H. Grauert and R. Remmert, {\it Theory of Stein spaces}, translated from the German by Alan Huckleberry. Reprint of the 1979 translation, 
Classics in Mathematics, Springer, Berlin, 2004; \href{https://mathscinet.ams.org/mathscinet/article?mr=2029201}{MR2029201}

\bibitem[GR84]{GR84} \bysame, {\it Coherent analytic sheaves}, Grundlehren der mathematischen Wissenschaften, 265, Springer, Berlin, 1984; \href{https://mathscinet.ams.org/mathscinet/article?mr=755331}{MR0755331}

\bibitem[Gro89]{Gro89} M. Gromov, Oka's principle for holomorphic sections of elliptic bundles, J. Amer. Math. Soc. {\bf 2} (1989), no.~4, 851--897; \href{https://mathscinet.ams.org/mathscinet/article?mr=1001851}{MR1001851}

\bibitem[H\"or90]{Hor90} L.~V. H\"ormander, {\it An introduction to complex analysis in several variables}, third edition, 
North-Holland Mathematical Library, 7, North-Holland, Amsterdam, 1990; \href{https://mathscinet.ams.org/mathscinet/article?mr=1045639}{MR1045639}

\bibitem[LL16]{LL16} R. L\"ark\"ang and F. L\'arusson, Extending holomorphic maps from Stein manifolds into affine toric varieties, Proc. Amer. Math. Soc. {\bf 144} (2016), no.~11, 4613--4626; \href{https://mathscinet.ams.org/mathscinet/article?mr=3544514}{MR3544514}

\bibitem[{\L}oj64]{Loj64} S. \L ojasiewicz, Triangulation of semi-analytic sets, Ann. Scuola Norm. Sup. Pisa Cl. Sci. (3) {\bf 18} (1964), 449--474; \href{https://mathscinet.ams.org/mathscinet/article?mr=173265}{MR0173265}

\bibitem[Oka39]{Oka39} K. Oka, Sur les fonctions analytiques de plusieurs variables. III---Deuxi\`eme probl\`eme de Cousin, J. Sci. Hiroshima Univ. Ser. A. Math. Phys. Chem. {\bf 9} (1939), 7--19; \href{https://mathscinet.ams.org/mathscinet/article?mr=4975236}{MR4975236}


\bibitem[Spa94]{Spa94} E.~H. Spanier, {\it Algebraic topology}, corrected reprint of the 1966 original, 
Springer, New York, 1994; \href{https://mathscinet.ams.org/mathscinet/article?mr=1325242}{MR1325242}

\end{thebibliography}
\end{document}